\documentclass[a4paper,11pt]{article}
\usepackage{amssymb,amsmath,bm,enumerate,theorem}
\usepackage[dvips]{graphicx}
\usepackage[all]{xy}

\renewcommand{\Im}{\text{Im}}
\newcommand{\Aut}{\mathrm{Aut}}

\newcommand{\Hom}{\mathrm{Hom}}
\newcommand{\Lin}{\text{Lin}}

\newcommand{\Der}{\text{Der}}

\newcommand{\R}{\mathcal{R}}
\newcommand{\D}{\mathcal{D}}
\newcommand{\M}{\mathcal{M}}
\newcommand{\C}{\mathcal{C}}
\newcommand{\g}{\mathfrak{g}}

\renewcommand{\l}{\mathfrak{l}}
\newcommand{\proof}{\textbf{Proof.}\;}
\newcommand{\deffont}{\it}
\renewcommand{\:}{\colon}
\theoremstyle{break}
\theorembodyfont{\normalfont}
\newtheorem{defi}{Definition}[section]
\newtheorem{lemm}[defi]{Lemma}

\newtheorem{prop}[defi]{Proposition}
\newtheorem{theo}[defi]{Theorem}
\newtheorem{remk}[defi]{Remark}
\newtheorem{fact}[defi]{Fact}

\newtheorem{conj}[defi]{Conjecture}

\newtheorem{sett}[defi]{Setting}
\newtheorem{nota}[defi]{Notation}

\title{A splitting of the local rigidity of Clifford-Klein forms of homogeneous spaces
of completely solvable Lie groups}
\author{Yoshinori Tanimura}
\begin{document}
\setlength{\abovedisplayskip}{5.5pt}
\setlength{\belowdisplayskip}{5.5pt}
\maketitle
\begin{abstract}
	In this article, we discuss the local rigidity of Clifford-Klein forms of homogeneous spaces of
	1-connected completely solvable Lie groups. In fact, we introduce a splitting of the local rigidity:
	vertical rigidity and horizontal rigidity. By using this splitting, we refine some
	existing results about the local rigidity and introduce a new approach to Baklouti's conjecture about
	the local rigidity.
\end{abstract}

\section{Introduction}
\subsection{Abstract}
For a Lie group $G$, a closed subgroup $H \subset G$ and a discrete subgroup $\Gamma \subset G$,
if the natural action $\Gamma \curvearrowright G/H$ is properly discontinuous and free,
then the double coset space $\Gamma \backslash G / H$ is called a Clifford-Klein form
(\cite{MR1016093}).
We want to know how the structure of Clifford-Klein form changes if $\Gamma$ is perturbed in $G$.
By formulating the classification space of Clifford-Klein forms,
the problem is boiled down to the topological structure of
$\D(\Gamma,G,H)$ and $\M(\Gamma,G,H)$ (Definition \ref{defi:moduli}).\par
In this article, we discuss local rigidity of Clifford-Klein forms (Definition \ref{defi:locally_rigid})
in the case of $G$ being nilpotent or completely solvable. In the reductive case, it is well-known
that there are some cases of structures of Clifford-Klein forms that satisfy local rigidity
(The Weil-Selberg-Kobayashi local rigidity theorem) or stronger rigidity (Mostow's rigidity theorem).
On the other hand, in the solvable case, structures of Clifford-Klein forms scarcely satisfy
local rigidity. In regards to this, the following conjecture is unsolved.

\begin{conj}[{\cite[Conjecture 3.1]{MR2863361}}] \label{conj:baklouti}
	Let $G$ be a 1-connected nilpotent Lie group, $H \subset G$ a closed subgroup and $\Gamma \subset G$
	a non-trivial discrete subgroup. Then no elements of $\R(\Gamma,G,H)$ are locally rigid.
\end{conj}
It is known that this conjecture is true in a special case.

\begin{fact}[{\cite[Theorem 3.5]{MR3043389}}] \label{fact:baklouti_conjecture}
	Let $G$ be a 1-connected nilpotent Lie group, $H \subset G$ a closed
	subgroup and $\Gamma \subset G$ a non-trivial discrete subgroup such that
	its syndetic hull is not characteristically nilpotent. Then no elements of $\R(\Gamma,G,H)$ are locally rigid.
\end{fact}
For example, if a 1-connected nilpotent Lie group $L$ is of dimension more than 1 and 2-step nilpotent,
abelian or of dimension lower than 7, then $L$ is not characteristically nilpotent.

In this article, we will refine Fact \ref{fact:baklouti_conjecture} and introduce a new approach to
Conjecture \ref{conj:baklouti} by splitting local rigidity to vertical rigidity and horizontal rigidity
(Definition \ref{defi:split_of_loc_rigid}).
In fact, Fact \ref{fact:baklouti_conjecture} will be made the following theorem.

\begin{theo} \label{theo:main}
	Let $G$ be a 1-connected nilpotent Lie group, $H \subset G$ a closed
	subgroup and $\Gamma \subset G$ a non-trivial discrete subgroup such that
	its syndetic hull is not characteristically nilpotent. Then no elements of $\R(\Gamma,G,H)$ are vertically rigid.
\end{theo}

\subsection{Background}
Some Clifford-Klein forms describe some geometrically important spaces.
So there are many researches being conducted about the classification of Clifford-Klein forms.
For example, the Teichm\"{u}ller theory, researched by Teichm\"{u}ller, Ahlfors, S. Kobayashi and so on,
is identified with the classification of Clifford-Klein forms which describe closed Riemann surfaces.
For the classification of Clifford-Klein forms which describe closed hyperbolic manifolds whose
dimension is no less than 3, Mostow's rigidity theorem is well-known. And for the classification of
Clifford-Klein forms which describe a certain kind of 3 dimensional Lorentzian manifold,
T. Kobayashi solved the Goldman conjecture.\par
For the classification of Clifford-Klein forms, the main interest is the cases of $G$ being simple,
semi-simple or reductive. This is because the cases of $G$ being simple, semi-simple or reductive
are often geometrically important. But if we recognize that the classification of
Clifford-Klein forms are the abstract problem for Lie groups, there is room to discuss
the case of $G$ not being reductive, especially of it being abelian, nilpotent or solvable.
From such reasons, some research works have recently been focused on the cases of $G$ being
abelian, nilpotent or solvable.

\section{Preliminaries}
Let $G$ be a Lie group, $H \subset G$ a closed subgroup and $\Gamma \subset G$ a discrete subgroup.
We define the classification spaces of properly discontinuous and free actions of
$\Gamma \curvearrowright G/H$.

\begin{defi} \label{defi:moduli}
	Let $\Hom(\Gamma,G)$ be the space of all homomorphisms of Lie groups with the pointwise-convergence
	topology and define the two actions $G \curvearrowright \Hom(\Gamma,G) \curvearrowleft \Aut(\Gamma)$
	by the adjoint action $G \curvearrowright G$ and the natural action
	$\Aut(\Gamma) \curvearrowright \Gamma$. We define these four spaces. \vspace{-5.5pt}
	\begin{itemize}
			\setlength{\leftskip}{-11pt}
			\setlength{\itemsep}{2.2pt}
		\item {\deffont The parameter space} (\cite[Remark 3]{MR1231232})\\
			$\R(\Gamma,G,H) := \left\{ \varphi \in \Hom(\Gamma,G) \Big|
			\begin{array}{l}
				\varphi \: \Gamma \hookrightarrow G \text{ is an embedding, and} \\
				\varphi(\Gamma) \curvearrowright G/H \text{ is properly discontinuous and free.}
			\end{array}
			\right\}$
		\item {\deffont The Chabauty space}\\
			$\C(\Gamma,G,H) := \R(\Gamma,G,H) / \Aut(\Gamma)$.
		\item {\deffont The deformation space} (\cite[Definition 5.3.1]{MR1852186}) \\
			$\D(\Gamma,G,H) := G \backslash \R(\Gamma,G,H)$.
		\item {\deffont The moduli space} (\cite[Definition 5.3.2]{MR1852186}) \\
			$\M(\Gamma,G,H) := G \backslash \R(\Gamma,G,H) / \Aut(\Gamma)$.
	\end{itemize}
	These spaces can be visually summarized as follows.
	\begin{align*}
		\xymatrix{
		\R(\Gamma,G,H) \ar[r]^{G \backslash} \ar[d]_{/ \Aut(\Gamma)}
		\ar@{}[dr]|{\text{\Large $\circlearrowleft$}}
		& \D(\Gamma,G,H) \ar[d]^{/ \Aut(\Gamma)} \\
		\C(\Gamma,G,H) \ar[r]_{G \backslash}
		& \M(\Gamma,G,H)
		}
	\end{align*}
\end{defi}

$\M(\Gamma,G,H)$ is the space of Clifford-Klein forms with a natural topology.
So for $\varphi \in \R(\Gamma,G,H)$, that $[\varphi]_{\M} \in \M(\Gamma,G,H)$ is isolated
means that the structure of Clifford-Klein form defined by $\varphi$ is invariant by perturbation
of $\varphi$. Moreover, since $\Aut(\Gamma)$ is discrete, we can guess the local structure of
$\D(\Gamma,G,H)$ is similar to that of $\M(\Gamma,G,H)$. So we define local rigidity as the following.

\begin{defi}[{\cite[Remark 3]{MR1231232}}] \label{defi:locally_rigid}
	An element $\varphi \in \R(\Gamma,G,H)$ is said to be {\deffont locally rigid} if
	$[\varphi]_{\D} \in \D(\Gamma,G,H)$ is isolated.
\end{defi}

\begin{remk}
	Later, we will discuss the case that $\Gamma$ is a connected Lie group. In this case, the topology of
	$\Hom(\Gamma,G)$ will be the compact-open topology.
\end{remk}

\section{A splitting of local rigidity and the going-through map}
\subsection{1-connected completely solvable Lie groups} \label{csl}
Let $G$ be a 1-connected completely solvable Lie group, $H$ a closed subgroup of $G$ and $\Gamma$ a
discrete subgroup of $G$.
As long as the main interest is in the deformation space, we can replace $\Gamma$ (resp. $H$) with a
closed connected subgroup $L$ (resp. $H$) of $G$ by considering the syndetic hull. Since $L$ is a
1-connected Lie group, $\R(L,G,H) \subset \Hom(L,G) \cong \Hom(\l,\g),\Aut(L) \cong \Aut(\l)$.

\begin{remk}
	The idea of using syndetic hulls appeared in \cite{MR2287673} for the first time. They discussed
	the case that $G$ is a specific nilpotent Lie group. A. Baklouti and I. K\'{e}dim discussed
	the general 1-connected nilpotent case in \cite{MR2548400}
	and the 1-connected completely solvable case in \cite{MR2609023}.
\end{remk}

In the following, we denote $\R(L,G,H)$ by $\R$ and respectively $\D,\C$ and $\M$ in the same way.
The $G$-equivalent principal bundle structure in the following proposition is essential.

\begin{prop} \label{prop:rispgb}
	The quotient map $\R \to \C$ is a principal $\Aut(\l)$-bundle.
\end{prop}
\proof We denote the set of all injective linear maps from $\l$ to $\g$ by $\Lin^{\circ}(\l,\g)$.
Then, since $\R \subset \Hom(\l,\g) \subset \Lin^\circ(\l,\g)$ and $\Aut(\l) \subset GL(\l)$,
$\R \to \C$ is a restriction of the frame bundle $\Lin^\circ(\l,\g) \to \Lin^\circ(\l,\g)/\Aut(\l)$.\par
\ \par
Later, we will introduce splitting of local rigidity by using
this bundle structure.

\subsection{On principal bundles} \label{pgb}
In this subsection, we discuss in the following condition for simplification of discussion.

\begin{sett} \label{sett:pgb}
	Let $X$ be a topological space, $G$ and $H$ topological groups
	and $X$ has a continuous left $G$-action $G \curvearrowright X$ and a continuous
	right $H$-action $X \curvearrowleft H$.
	We assume the actions $G \curvearrowright X$ and $X \curvearrowleft H$ are commutative and
	the quotient map $X \to X/H$ a principal $H$-bundle.
	And let ${}_G \pi,\pi_H,{}_G \varpi$ and $\varpi_H$ in the following diagram be the quotient maps.
	\begin{align*}
		\xymatrix{
			X \ar[r]^{ {}_G \pi } \ar[d]_{\pi_H}
			\ar@{}[dr]|{\text{\Large $\circlearrowleft$}}
			& G \backslash X \ar[d]^{\varpi_H} \\
			X/H \ar[r]_{ {}_G \varpi }
			& G \backslash X / H
		}
	\end{align*}
\end{sett}

\subsubsection{Splitting of local rigidity}
At first, we will split the local rigidity into two conditions.

\begin{defi} \label{defi:split_of_loc_rigid}
	Assume Setting \ref{sett:pgb} and take $x \in X$. \vspace{-3.3pt}
	\begin{itemize}
			\setlength{\leftskip}{-11pt}
			\setlength{\itemsep}{2.2pt}
		\item $x$ is {\deffont locally rigid} if $Gx \in G \backslash X$ is isolated.
		\item $x$ is {\deffont vertically rigid} if the $H$-orbit of $Gx \in G \backslash X$ is discrete.
		\item $x$ is {\deffont horizontally rigid} if $GxH \in G \backslash X / H$ is isolated.
	\end{itemize}
\end{defi}

\begin{remk}
	Assume Setting \ref{sett:pgb}. For $x \in X$, the followings are equivalent. \vspace{-3.3pt}
	\begin{enumerate}[{(}1{)}]
			\setlength{\itemsep}{2.2pt}
		\item $x$ is locally rigid.
		\item $x$ is vertically and horizontally rigid.
	\end{enumerate}
\end{remk}

\subsubsection{The going-through map}
Next we introduce the going-through map. The going-through map is useful to discuss vertical rigidity.
More precisely, we can describe a necessary condition of vertical rigidity by using the going-through map.

\begin{nota} \label{nota:pgb}
	Assume Setting \ref{sett:pgb}. Take a point $x \in X$. \vspace{-3.3pt}
	\begin{itemize}
			\setlength{\leftskip}{-11pt}
			\setlength{\itemsep}{2.2pt}
		\item $G_{xH}$ (resp. ${}_{Gx} H$) denotes the isotropy group of the action
			$G \curvearrowright X/H$ (resp. $G \backslash X \curvearrowleft H$) at $xH$ (resp. $Gx$).
		\item We define two continuous maps ${}_G \theta_x \: G \to X$ and $\theta_{x,H} \: H \to X$ by
			the followings:
			\begin{align*}
				{}_G \theta_x(g) := gx \; (g \in G) ,\; \theta_{x,H}(h) := xh \; (h \in H)
			\end{align*}
			Since the quotient map $H \to X/H$ is a principal $H$-bundle,
			$\theta_{x,H} \: H \to x \cdot H$ is a homeomorphism.
		\item Let $\bar{\theta}_{x,H} \: {}_{Gx} H \backslash H \to G \backslash X$ be the map
			induced by ${}_G \pi \circ \theta_{x,H} \: H \to G \backslash X$. Then the map
			$\bar{\theta}_{x,H} \: {}_{Gx} H \backslash H \to (Gx) \cdot H$ is bijective and continuous
			(however not necessarilly homeomorphism (See Subsection \ref{subs:fiber})).
	\end{itemize}
\end{nota}

\begin{lemm} \label{lemm:going_through}
	For $x \in X$, there exists a unique surjective continuous homomorphism
	$\alpha_x \: G_{xH} \to {}_{Gx} H$ such that $gx = x \alpha_x(g) \; (g \in G_{xH})$.
\end{lemm}
\proof Since the action $X \curvearrowleft H$ is free, $\alpha_x$ is unique.
And the composition of ${}_G \theta_x|_{G_{xH}} \: G_{xH} \to x \cdot H$ and
$\theta_{x,H}^{-1} \: x \cdot H \to H$ satisfies the condition of $\alpha_x$.

\begin{defi} \label{defi:going_through}
	We call the map $\alpha_x$ defined in Lemma \ref{lemm:going_through}
	{\deffont the going-through map of $x$}.
\end{defi}

In the case of the deformation space, the going-through map is a same map defined in
\cite[Equation 3.2]{MR3043389}. By using the going-through map, we obtain the following observation.
We denote the identity component of $H$ by $H_0$.

\begin{prop} \label{prop:vertical}
	Assume Setting \ref{sett:pgb}.
	If there exists a locally rigid point $x \in X$, then $H_0$ is a subquotient group of $G$.
\end{prop}
\proof
Since $\alpha_x \: G_{xH} \to {}_{Gx} H$ is surjective, it is enough to show $H_0 \subset {}_{Gx} H$.
Since the natural map ${}_{Gx} H \backslash H \to (Gx) \cdot H \subset G \backslash X$ is
bijective continuous and $(Gx) \cdot H \subset G \backslash X$ is discrete, ${}_{Gx} H \backslash H$ is
also discrete.

\subsection{On the deformation spaces}
Now we will make use of the notation of Subsection \ref{csl}. By Lemma \ref{prop:rispgb},
$(G,\R,\Aut(\l))$ satisfies the condition of $(G,X,H)$ in Setting \ref{sett:pgb}.
So we can define local rigidity, horizontal rigidity, vertical rigidity and the going-through maps for
elements of $\R$. And the local rigidity defined here is equivalent to the local rigidity defined in
Definition \ref{defi:locally_rigid}. By Proposition \ref{prop:vertical}, we get the following.

\begin{theo} \label{theo:main_precise}
	Let $G$ be a completely solvable (resp. nilpotent) Lie group, $H \subset G$ a closed subgroup,
	$\Gamma \subset G$ a discrete subgroup and $L$ a syndetic hull of $\Gamma$.
	Then if the derivation algebra of $\l$ is NOT completely solvable (resp. nilpotent),
	then no elements of $\R(\Gamma,G,H)$ are vertically rigid.
\end{theo}

A Lie algebra is called {\deffont characteristically nilpotent} if the Lie algebra of its derivations is
nilpotent. So Theorem \ref{theo:main_precise} is a precise description of Theorem \ref{theo:main}.
The details of characteristically nilpotent Lie algebras are written in \cite{MR1880104}.

\begin{remk}
	For $\varphi \in \R$, the isotropy group $G_{\varphi \Aut(\l)} \subset G$
	at $\varphi \Aut(\l) \in \C$ is the normalizer of the Lie subalgebra
	$\varphi(\l) \subset \g$. In particular, the condition $G_{\varphi \Aut(\l)} = G$ is equivalent to
	that $\varphi(\l) \subset \g$ is an ideal.
\end{remk}

\subsection{A sufficient condision that $\bar{\theta}_{x,H}$ is homeomorphism} \label{subs:fiber}
In Setting \ref{sett:pgb}, the map of $\bar{\theta}_{x,H}$ in Notation \ref{nota:pgb} is NOT
homeomorphic in generally. We want a sufficient condition that $\bar{\theta}_{x,H}$ is homeomorphic.

\begin{prop} \label{prop:nescon}
	Assume Setting \ref{sett:pgb}. Let $A \subset G \backslash X/H$ be a subset,
	$P := {}_G \varpi^{-1}(A) \subset X/H$, 
	$s \: P \to X$ a section of $\pi_H \: X \to X/H$ and $H' \subset H$
	a closed subgroup. Assume the following conditions.
	\begin{enumerate}[{(}1{)}]
			\setlength{\itemsep}{2.2pt}
		\item For all $p \in P$ and $g \in G$,
			there exists $g_0 \in G$ such that $s(gp) = g_0s(p)$.
		\item For all $p \in P$, ${}_{Gs(p)} H = H'$.
	\end{enumerate}
	Then there exists the unique homeomorphism $\Phi_s \: A \times (H' \backslash H)
	\stackrel{\cong}{\to} \varpi_H^{-1}(A)$
	such that $\Phi_s({}_G \varpi(p),H'h) = \bar{\theta}_{s(p),H}(H'h) \;
	(p \in P,\; h \in H)$.
\end{prop}

\begin{remk}
	If $P$ consists one point, the conditions (1) and (2) are true.
\end{remk}
\underline{\textbf{Proof of Proposition \ref{prop:nescon}}}\ \ 
Assume the following claims for all element $(p,h) \in P \times H$.
\begin{enumerate}
		\setlength{\leftskip}{-11pt}
		\setlength{\itemsep}{2.2pt}
	\item For all $g \in G$, there exists $h' \in H'$ such that $gs(p)h = s(gp)h'h$.
	\item For all $g \in G, h' \in H'$, there exists $g' \in G$ such that $g's(p)h = s(gp)h'h$.
\end{enumerate}
Then the following equivalent relations are the same.
\begin{itemize}
		\setlength{\leftskip}{-11pt}
		\setlength{\itemsep}{2.2pt}
	\item The relation defined by two actions $G \curvearrowright P$ and $H' \curvearrowright H$.
	\item The relation defined by the action $G \curvearrowright P \times H$ induced by
		the identification $X|_P \cong P \times H$ with respect to $s$.
\end{itemize}
The proposition follows from this observation immediately.\par
We will prove above two clames. By the assumption (1), there exists a map $\beta \: G \cdot p \to G$
such that $s(gp) = \beta(gp)s(p) \; (g \in G)$.\par
At first, we prove the claim 1. Let $g \in G$. Since $g \beta(gp)^{-1} \in G_{gp} = G_{s(gp)H}$,
\begin{align*}
	gs(p)h = g \beta(gp)^{-1} s(gp) h = s(gp) \cdot \alpha_{s(gp)} (g \beta(gp)^{-1}) h.
\end{align*}
Since $\alpha_{s(gp)}(g \beta(gp)^{-1}) \in {}_{Gs(gp)} H = H'$, the claim 1 holds.\par
Next we prove claim 2. Let $g \in G$ and $h' \in H$. By the surjectivity of the going-through map
$\alpha_{s(p)} \: G_{s(p)H} \to {}_{Gs(p)} H = H'$, there exists $g_0 \in G_{s(p)H}$ such that
$h' = \alpha_{s(p)}(g_0)$. So
\begin{align*}
	s(gp)h'h = \beta(gp) s(p) \alpha_{s(p)}(g_0) h = \beta(gp) g_0 s(p) h.
\end{align*}
This implies the claim 2.

\subsection{Examples of vertical rigidity and horizontal rigidity}
Ali Baklouti et al. conjectured there are no locally rigid Clifford-Klein forms in the nilpotent case
(Conjecture \ref{conj:baklouti}). However, there are some vertically rigid ones and horizontally rigid ones.

\subsubsection{An example of vertical rigidity}
At first, we construct a vertically rigid Clifford-Klein forms.

\begin{prop} \label{prop:example_vertical}
	Let $\l$ be a characteristically nilpotent Lie algebra of dimension more than 1, $\g := \Der(\l) \ltimes \l$,
	$L$ (resp. $G$) the 1-connected Lie group accociated by $\l$ (resp. $\g$) and $H := \{e\}$.
	Then the natural inclusion $\varphi \: \l \hookrightarrow \g$ is vertically rigid in $\R$.
\end{prop}

To prove this proposition, we use the following fact about characteristic nilpotency.

\begin{fact}[{\cite[Theorem 1]{MR0114841}}]
	For a Lie algebra $\l$ of dimension more than 1, the followings are equivalent.
	\begin{enumerate}[{(}1{)}]
			\setlength{\itemsep}{2.2pt}
		\item $\Der(\l)$ is nilpotent (i.e. $\l$ is characteristically nilpotent).
		\item $\Der(\l)$ consists nilpotent elements.
		\item $\Der(\l) \ltimes \l$ is nilpotent.
	\end{enumerate}
	In particular, if $\l$ is characteristically nilpotent, $\Aut(\l)_0$ is 1-connected.
\end{fact}
\underline{\textbf{Proof of Proposition \ref{prop:example_vertical}}}\ \ 
Since $\varphi(\l) \subset \g$ is an ideal, $[\varphi]_{\C} \in \C$ is a fixed point of the
action $G \curvearrowright \C$. So we obtain a homeomorphism $(\Im \alpha_\varphi) \backslash \Aut(\l)
\cong [\varphi]_{\D} \cdot \Aut(\l)$ by Proposition \ref{prop:nescon}.
Under the ideintification $\Aut(\l) \cong \Aut(L)$, since $\alpha_{\varphi}(T,0) = T$ for all
$(T,0) \in G = \Aut(\l)_0 \ltimes L$, $\Im \alpha_\varphi = \Aut(\l)_0$. So
$[\varphi]_{\D} \cdot \Aut(\l) \cong \Aut(\l)_0 \backslash \Aut(\l)$. It means $\varphi \in \R$ is
vertically rigid.

\subsubsection{Examples of horizontal rigidity}
Next we construct a horizontally rigid Clifford-Klein form. In the following, we denote $\R(L,G,\{e\})$ by
$\R(L,G)$.

\begin{prop}
	Let $G$ be a 1-connected completely solvable Lie group. Then
	$\R(G,G) = \Aut(\g), \D(G,G) = \textrm{Out}(\g)$ and each $\C(G,G)$ and $\M(G,G)$ consists one point.
	In particular, all elements of $\R(G,G)$ are horizontally rigid.
\end{prop}

Since the proof of this proposition is easy, we omit it. In this example, we consider compact
Clifford-Klein forms. On the other hand, there is an example of horizontal rigidity
of non-compact Clifford-Klein forms. Let $n$ be a positive integer, define
a linear map $\sigma_n \: \mathbb{R}^n \to \mathbb{R}^n$ by
$\sigma_n(x_0,\dots,x_{n-1}) = (0,x_0,\dots,x_{n-2})$ and let
$\l_n := \mathbb{R} \sigma_n \ltimes \mathbb{R}^n$. Then $\l_n$ is called the $n+1$ dimensional
ladder Lie algebra. Let $L_n$ be the 1-connected Lie group associated to $\l_n$.
Then the following proposition can be shown immediately.

\begin{prop}
	If $n > 3$, $(0,\mathbb{R}^n) \subset \l_n$ is only the $n$-dimensional abelian
	Lie subalgebra.
	In particular, $\C(\mathbb{R}^n,L_n)$ consists one point and
	all elements of $\R(\mathbb{R}^n,L_n)$ are horizontally rigid.
\end{prop}

\bibliographystyle{alpha}
\bibliography{references}

\begin{thebibliography}{Kob01}

\bibitem[AC01]{MR1880104}
Jos{\'e}~M. Ancochea and Rutwig Campoamor.
\newblock Characteristically nilpotent {L}ie algebras: a survey.
\newblock {\em Extracta Math.}, 16(2):153--210, 2001.

\bibitem[Bak11]{MR2863361}
Ali Baklouti.
\newblock On discontinuous subgroups acting on solvable homogeneous spaces.
\newblock {\em Proc. Japan Acad. Ser. A Math. Sci.}, 87(9):173--177, 2011.

\bibitem[BK09]{MR2548400}
Ali Baklouti and Imed K{\'e}dim.
\newblock On the deformation space of {C}lifford-{K}lein forms of some
  exponential homogeneous spaces.
\newblock {\em Internat. J. Math.}, 20(7):817--839, 2009.

\bibitem[BK10]{MR2609023}
Ali Baklouti and Imed K{\'e}dim.
\newblock On non-abelian discontinuous subgroups acting on exponential solvable
  homogeneous spaces.
\newblock {\em Int. Math. Res. Not. IMRN}, (7):1315--1345, 2010.

\bibitem[BK13]{MR3043389}
Ali Baklouti and Imed Kedim.
\newblock On the local rigidity of discontinuous groups for exponential
  solvable {L}ie groups.
\newblock {\em Adv. Pure Appl. Math.}, 4(1):3--20, 2013.

\bibitem[KN06]{MR2287673}
Toshiyuki Kobayashi and Salma Nasrin.
\newblock Deformation of properly discontinuous actions of $\mathbb{Z}^k$ on
  $\mathbb{R}^{k+1}$.
\newblock {\em Internat. J. Math.}, 17(10):1175--1193, 2006.

\bibitem[Kob89]{MR1016093}
Toshiyuki Kobayashi.
\newblock Proper action on a homogeneous space of reductive type.
\newblock {\em Math. Ann.}, 285(2):249--263, 1989.

\bibitem[Kob93]{MR1231232}
Toshiyuki Kobayashi.
\newblock On discontinuous groups acting on homogeneous spaces with noncompact
  isotropy subgroups.
\newblock {\em J. Geom. Phys.}, 12(2):133--144, 1993.

\bibitem[Kob01]{MR1852186}
Toshiyuki Kobayashi.
\newblock Discontinuous groups for non-{R}iemannian homogeneous spaces.
\newblock In {\em Mathematics unlimited---2001 and beyond}, pages 723--747.
  Springer, Berlin, 2001.

\bibitem[LT59]{MR0114841}
G.~Leger and S.~T{\^o}g{\^o}.
\newblock Characteristically nilpotent {L}ie algebras.
\newblock {\em Duke Math. J.}, 26:623--628, 1959.

\end{thebibliography}

\ \par
Graduate School of Mathematical Sciences, The University of Tokyo\par
This work was supported by the Program for Leading Graduate Schools, MEXT, Japan.

\end{document}